\documentstyle[amscd,amssymb,verbatim,12pt]{amsart}

\theoremstyle{plain}
\newtheorem{thm}{Theorem}
\newtheorem{prop}[thm]{Proposition}
\newtheorem{lem}[thm]{Lemma}
\newtheorem{cor}[thm]{Corollary}

\theoremstyle{definition}

\newenvironment{pfm}{\noindent{\em Proof of Theorem \ref{main}.}}{\qed}

\newcommand{\msp}{(\Omega,\Sigma,\mu)}
\newcommand{\mspp}{(\Omega ',\Sigma ',\mu ')}
\newcommand{\wLp}{L^{p,\infty}}
\newcommand{\wlp}{\ell^{p,\infty}}
\newcommand{\wLt}{L^{2,\infty}}
\newcommand{\wlt}{\ell^{2,\infty}}
\newcommand{\Mp}{M^{p,\infty}}
\newcommand{\mpi}{m^{p,\infty}}
\newcommand{\Mt}{M^{2,\infty}}
\newcommand{\mti}{m^{2,\infty}}
\newcommand{\al}{\alpha}

\newcommand{\Om}{\Omega}
\newcommand{\si}{\sigma}
\newcommand{\Si}{\Sigma}
\newcommand{\ep}{\epsilon}
\newcommand{\G}{\Gamma}
\newcommand{\g}{\gamma}
\newcommand{\dsum}{\oplus_{\al\in A}}
\newcommand{\two}{\{-1,1\}}
\newcommand{\ap}{\aleph}

\newcommand{\spn}{\operatorname{span}}
\newcommand{\supp}{\operatorname{supp}}
\newcommand{\sgn}{\operatorname{sgn}}
\newcommand{\la}{\langle}
\newcommand{\ra}{\rangle}

\newcommand{\N}{{\Bbb N}}
\newcommand{\R}{{\Bbb R}}
\newcommand{\bs}{\backslash}

\newcommand{\bl}{\bigl}
\newcommand{\br}{\bigr}
\newcommand{\bc}{\bigcup}

\begin{document}


\title{Purely non-atomic weak $L^p$ spaces}
\author{Denny H.\ Leung}
\address{Department of Mathematics\\
         National University of Singapore\\
         Singapore 119260}
\email{matlhh@@leonis.nus.sg}



\maketitle

\begin{abstract}
Let $\msp$ be a purely non-atomic measure space, and let $1 < p <
\infty$. If $\wLp\msp$ is isomorphic, as a Banach space, to $\wLp\mspp$
for some purely atomic measure space $\mspp$, then there is a
measurable partition $\Omega = \Omega_1\cup\Omega_2$ such that
$(\Omega_1,\Sigma\cap\Omega_1,\mu_{|\Sigma\cap\Omega_1})$ is countably
generated and $\sigma$-finite, and that $\mu(\sigma) = 0$ or $\infty$
for every measurable $\sigma \subseteq \Omega_2$. In particular,
$\wLp\msp$ is isomorphic to $\ell^{p,\infty}$. 
\end{abstract}



\section{Introduction}
In \cite{L1}, the author proved that the spaces $\wLp [0,1]$ and $\wLp
[0,\infty)$ are both isomorphic to the atomic space $\wlp$.
Subsequently, it was observed that if $\msp$ is countably generated and
$\sigma$-finite, then $\wLp\msp$ is isomorphic to either $\wlp$ or
$\ell^\infty$ \cite[Theorem 7]{L2}. 
In this paper, we show that the isomorphism of atomic
and non-atomic weak $L^p$ spaces does not hold beyond the countably
generated, $\sigma$-finite situation. 

Before giving the precise statement of the main theorem, let us agree
on some terminology.  Throughout this paper, every measure space under
discussion is assumed to be {\em non-trivial}\/ in the sense that it
contains a measurable subset of finite non-zero measure. A measurable
subset $\sigma$ of a measure space $\msp$ is an {\em atom}\/ if
$\mu(\sigma) > 0$, and either
$\mu(\sigma') = 0$ or $\mu(\sigma\backslash\sigma') = 0$ for each
measurable subset $\sigma'$ of $\sigma$. A {\em purely non-atomic}\/
measure space is one which contains no atoms.
We say that a collection $S$ of measurable sets {\em generates}\/ a
measure space $\msp$ if $\Sigma$ is the smallest $\sigma$-algebra
containing $S$ as well as the $\mu$-null sets.
A measure space $\msp$ is {\em purely atomic}\/ if it is generated by
the collection of all of its atoms; it is 
{\em countably generated}\/ if there is a sequence $(\sigma_n)$ in
$\Sigma$ which generates $\msp$.
For any measure space $\msp$, and $1 < p < \infty$, the {\em weak $L^p$
space}\/ $\wLp\msp$ is the space of all (equivalence classes of)
$\Sigma$-measurable functions
$f$ such that 
\[ \|f\| = \sup_{c>0}c(\mu\{|f| > c\})^{\frac{1}{p}} < \infty . \]
It is well known that $\|\cdot\|$ is equivalent to a norm under which
$\wLp\msp$ is a Banach space.  However, since we are only concerned
with isomorphic questions, we will employ the quasi-norm $\|\cdot\|$
exclusively in our computations.
The aim of this paper is to prove the following theorem.

\begin{thm}\label{main}
Let $\msp$ be a purely non-atomic
measure space, and let $1 < p < \infty$. The following
statements are equivalent.
\begin{enumerate}
\item $\wLp\msp$ is isomorphic to $\wLp\mspp$ for some purely atomic
measure space $\mspp$. 
\item $\wLp\msp$ is isomorphic to a subspace of 
$\wLp\mspp$ for some purely atomic measure space $\mspp$. 
\item There is a measurable partition $\Omega = \Omega_1 \cup \Omega_2$
such that $(\Omega_1,\Sigma\cap\Omega_1,\mu_{|\Sigma\cap\Omega_1})$ is
countably
generated and $\sigma$-finite, and that $\mu(\sigma) = 0$ or $\infty$
for every measurable $\sigma \subseteq \Omega_2$.
\item $\wLp\msp$ is isomorphic to $\wlp$.
\end{enumerate}
\end{thm}

It is interesting to note that with regard to (2), the weak $L^p$
spaces behave in a way that is ``in between'' the behavior of the $L^p$
spaces, $1 \leq p < \infty$, and $L^\infty$. Indeed, if $\msp$ is purely
non-atomic,  then $L^p\msp$ can never be embedded into an atomic $L^p$
space ($1 \leq p < \infty,\,\, p \neq 2$). On the other hand, along
with all Banach spaces, $L^\infty\msp$ is isomorphic to a subspace of
$\ell^\infty(J)$ for a sufficiently large index set $J$.  

The other notation follows mainly that of \cite{LT,LT2}. Banach spaces
$E$ and $F$ are said to be {\em isomorphic}\/ if they are linearly
homeomorphic; $E$ {\em embeds}\/ into $F$ if it is isomorphic to a
subspace of $F$. If $I$ is an arbitrary index set, and $(x_i)_{i\in
I}$, $(y_i)_{i\in I}$ are indexed collections of elements in possibly
different Banach spaces, we say that they are {\em equivalent}\/ if
there is a constant $0 < K < \infty$ such that 
\[ K^{-1}\bl\|\sum a_ix_i\br\| \leq \bl\|\sum a_iy_i\br\| 
\leq K\bl\|\sum a_ix_i\br\| \]
for every collection $(a_i)_{i\in I}$ of scalars with finitely many
non-zero terms.
We will also have occasion to use terms and notation
concerning vector lattices, for which the references are
$\cite{LT2,S}$. In particular, two elements $a, b$ of a vector lattice
are said to be {\em disjoint}\/ if $|a|\wedge|b| = 0$. A Banach lattice
$E$ satisfies an {\em upper $p$-estimate}\/ if there is a constant $M <
\infty$ such that 
\[ \bl\|\sum^n_{i=1}x_i\br\| 
\leq M\bl(\sum^n_{i=1}\|x_i\|^p\br)^{\frac{1}{p}} \]
whenever $(x_i)^n_{i=1}$ is a pairwise disjoint sequence in $E$.  It is
trivial to check that every $\wLp\msp$ satisfies an upper $p$-estimate
with constant $1$. Finally, if $A$ is an arbitrary set, we write 
${\cal P}(A)$ for the power set of $A$, and $A^c$ for its complement
(with respect to some universal set). 
 
\section{Proof of the Main Theorem}

Let us set the notation for the two types of measure spaces which will
command a large part of our attention. By $\{-1,1\}$, we will mean the
two-point measure space, each point of which is assigned a mass of
$\frac{1}{2}$. If $I$ is an arbitrary index set, $\{-1,1\}^I$ is the
product measure space of $I$ copies of $\{-1,1\}$. Now let
$((\Om_\al,\Sigma_\al,\mu_\al))_{\al\in A}$ be a collection of pairwise
disjoint measure spaces. We define the measurable space $(\Om,\Sigma)$
to be the set $\cup_{\al\in A}\Om_\al$, endowed with the smallest
$\sigma$-algebra $\Sigma$ generated by $\cup_{\al\in A}\Sigma_\al$. For
any $\sigma \in \Sigma$, define 
\[ \mu(\sigma) = \sum_{\al\in A}\mu_\al(\sigma\cap{\Om_\al}) .\]
The measure space $(\Om,\Sigma,\mu)$ is denoted by $\oplus_{\al\in
A}(\Om_\al,\Sigma_\al,\mu_\al)$. Of particular interest will be
$\dsum J_\al$, where
each $J_\al$ is a copy of the measure space $[0,1]$ with the Lebesgue
measure.

\begin{thm}\label{short}
If $I$ is an uncountable index set, $\wLp(\{-1,1\}^I)$ does not
embed into $\wLp\msp$ for any purely atomic measure
space $\msp$.
\end{thm}

\begin{thm}\label{long}
Let $A$ be an arbitrary index set. For every $\al \in A$, let
$J_\al$ be a copy of the measure space $[0,1]$.
If $\wLp(\dsum J_\al)$ embeds into $\wLp\msp$ for
some purely atomic measure space $\msp$, then the set $A$ is countable.
\end{thm}

The proofs of the crucial Theorems \ref{short} and \ref{long} will be
the subject of the subsequent sections. 
To apply these theorems to the proof of the
main theorem (Theorem \ref{main}) requires the use of certain known
facts, which we now recall. 
Let $\msp$ and $\mspp$ be measure spaces. Denote by $\Theta_\mu$ and
$\Theta_{\mu'}$ the $\mu$- and $\mu'$-null sets respectively. Then
$\mu$ induces a function on the $\sigma$-complete Boolean algebra
$\Sigma/\Theta_{\mu}$, which we denote again by $\mu$. Similarly for
$\mu'$. We say that the measure spaces $\msp$ and $\mspp$ are
{\em isomorphic}\/ if there exists a Boolean algebra isomorphism $\Phi
: \Sigma/\Theta_{\mu} \to \Sigma'/\Theta_{\mu'}$ such that $\mu =
\mu'\circ\Phi$. For notions and results regarding measure algebras, we
refer to \cite[\S 14]{La}. 
The next fact, which can be found in
\cite{Lo}, follows easily from the observation that the set of
functions $f \in \wLp\msp$ of the form $f = \bigvee a_n\chi_{A_n}$,
where $(a_n) \subseteq \R$, and $(A_n)$ is a pairwise disjoint sequence
in $\Sigma$, is dense in $\wLp\msp$.

\begin{thm}\label{iso}
If $\msp$ and $\mspp$ are isomorphic measure spaces, 
then the Banach spaces
$\wLp\msp$ and $\wLp\mspp$ are isometrically isomorphic.
\end{thm}

The next theorem is stated in the form in which we will use it.  It is
a consequence of Maharam's theorem on the classification of measure
algebras; see \cite[Theorems 14.7 and 14.8]{La}. If $\msp$ is a measure
space, and $c$ is a positive number, we let $c\mu$ be the measure given by
$(c\mu)(\sigma) = c \mu(\sigma)$ for all $\sigma \in \Sigma$. Clearly,
the map $f \mapsto c^{-1/p}f$ is an isometric
isomorphism from $\wLp\msp$ onto $\wLp(\Om,\Si,c\mu)$. 

\begin{thm}[Maharam]\label{maharam}
Let $\msp$ be a purely non-atomic, finite measure space which is not
countably generated. Then there is a measurable subset $\Om'$ of $\Om$
and an uncountable index set $I$ such that 
$\mspp$ is isomorphic to $\{-1,1\}^I$, 
where $\Sigma' = \Sigma\cap\Om'$, and $\mu' =
(\mu(\Om'))^{-1}\mu_{|\Sigma'}$. 
Consequently, $\wLp\msp$ has a subspace isometrically isomorphic to 
$\wLp(\{-1,1\}^I)$.
\end{thm}  

For the following proof, recall that a Banach space $E$ satisfies the
{\em Dunford-Pettis property}\/ if $\la x_n,x'_n\ra \to 0$ whenever
$(x_n)$ and $(x'_n)$ are weakly null sequences in $E$ and $E'$
respectively. It is well known that $\ell^\infty$ satisfies the
Dunford-Pettis property; see, e.g., \cite[\S II.9]{S}.\\

\begin{pfm}
Suppose (3) holds.  Then $\wLp\msp$ is clearly isometrically
isomorphic to
$\wLp(\Omega_1,\Sigma\cap\Omega_1,\mu_{|\Sigma\cap\Omega_1})$.  By
\cite[Theorem 7]{L2}, $\wLp\msp$ is isomorphic to either $\wlp$ or
$\ell^\infty$. However, since $\msp$ is purely non-atomic, we can
easily verify that $\wLp\msp$ fails the Dunford-Pettis property. (Use
Rademacher-like functions.) Hence it cannot be isomorphic to $\ell^\infty$.
The implications (4) $\Rightarrow$ (1) $\Rightarrow$ (2) are trivial.
Therefore, it remains to prove that (2) $\Rightarrow$ (3).  Using
Zorn's Lemma, obtain a (possibly empty) collection of measurable subsets
$(\Om_\al)_{\al\in A}$ of $\Om$ which is maximal with respect to the
following conditions : 
$\mu(\Om_\al\cap\Om_\beta) = 0$ if $\al \neq \beta$;
$\mu(\Om_\al) = 1$ for all $\al\in A$.
For each $\al \in A$, let $J_\al$ be a copy of the measure space
$[0,1]$. Then $J_\al$ is isomorphic to a measure subalgebra of
$(\Om_\al,\Sigma\cap\Om_\al,\mu_{|\Sigma\cap\Om_\al})$.  
It follows that $\oplus_{\al\in A}J_\al$ is isomorphic to a measure
subalgebra of $\msp$. Theorem \ref{iso} implies that $\wLp(\dsum
J_\al)$ is isometrically isomorphic 
to a subspace of $\wLp\msp$, and hence, by the
assumption (2), isomorphic to a subspace of an atomic weak $L^p$
space. According to Theorem \ref{long}, $A$ must be a countable set.
By the maximality of $(\Om_\al)_{\al\in A}$\,, 
\begin{multline*} 
m \equiv \sup\{\mu(\si) :
\si\ \text{is a measurable subset of}\\
\text{$\Om\backslash\!\cup_{\al\in A}\Om_\al$ of finite measure}\} 
\leq 1. 
\end{multline*}
It is easily seen that the supremum is attained, say, at $\Om_0$.
Define $\Om_1 = \Om_0\cup(\cup_{\al\in A}\Om_\al)$.  Since $A$ is
countable, $\Om_1 \in \Si$. If $\Om_\al$ is not countably generated for
some $\al\in A$, then Theorem \ref{maharam} produces
an uncountable index set $I$ such that $\wLp(\{-1,1\}^I)$ is 
isometrically isomorphic
to a subspace of $\wLp(\Om_\al)$, and thus isomorphic to a subspace of
an atomic weak $L^p$ space. This violates Theorem \ref{short}.
Similarly, we see that $\Om_0$ is countably generated. Therefore,
$\Om_1$ is countably generated; it is clearly $\si$-finite. 
If $\si$ is a measurable subset of $\Om_2 = \Om\backslash\Om_1$, 
and $0 < \mu(\si) < \infty$, then  $m < \mu(\Om_0\cup\si) < \infty$, 
contrary to
the choice of $\Om_0$. Hence $\mu(\si) = 0$ or $\infty$.
\end{pfm}
 
\section{The space $\wLp(\{-1,1\}^I)$}

Let $\Gamma$ be an arbitrary set, and let $w : \Gamma \to (0,\infty)$ be
a weight function. We can define a measure $\mu$ on ${\cal P}(\Gamma)$
by $\mu(\si) = \sum_{\gamma\in\,\si}w(\gamma)$
for all $\si \subseteq \Gamma$. The resulting weak $L^p$ space 
$\wLp(\Gamma,{\cal P}(\Gamma),\mu)$ will be denoted by $\wlp(\G,w)$,
or simply $\wlp(\G)$ if $w$ is identically $1$. It
is easy to see that if $\msp$ is purely atomic, then $\wLp\msp$ is 
isometrically isomorphic 
to $\wlp(\G,w)$ for some $(\G,w)$. If $\msp$ is a measure space,
and $1 < p < \infty$, let $\Mp\msp$ be the closed subspace of $\wLp\msp$
generated by 
the functions $\chi_\si$, where $\si$ is a measurable set of finite measure. 
The corresponding subspace of $\wlp(\G,w)$ is 
denoted by $\mpi(\G,w)$. The proof of Theorem \ref{short} for the case
$ p \neq 2$ is rather easy and is contained in Theorem \ref{l2I}. 
For the reader's convenience, we recall the following 
disjunctification
result \cite[Proposition 10]{L2}.

\begin{prop}\label{disj}
Let $w$ be a weight function on a set $\G$. Assume that $A$ and $B$ are
subsets of $\wlp(\G,w)$
such that $|A| > \max\{|B|,\aleph_0\}$. Suppose also that there are
constants $K < \infty$, $r > 1$ such that 
\begin{equation} \label{dom}
 \bl\|\sum_{x\in F}\ep_xx\br\| \leq K|F|^{\frac{1}{r}} 
\end{equation}
for all finite subsets $F$ of $A$, and all $\ep_x = \pm 1$. Then there
exists $C \subseteq A$, $|C| = |A|$, such that the elements of $C$ are 
pairwise disjoint, and $|b|\wedge|c| = 0$ whenever $b \in B$, $c \in C$.
\end{prop} 

\begin{pf}
First we show that if $\G'$ is a subset of $\G$ such that $|\G'| <
|A|$, then there exists $A' \subseteq A$, $|A'| = |A|$, such that
$x\chi_{\G'} = 0$ for all $x \in A'$.
Indeed, let $A' = \{x \in A : x\chi_{\G'} = 0\}$.  
For each $x \in A\bs A'$,
there is a $\g \in \G'$ such that $x(\g) \neq 0$.  Define a choice
function $f : A\bs A' \to \G'$ such that $x(f(x)) \neq 0$ for all $x
\in A \bs A'$.  If $|A'| < |A|$, then $|A\bs A'| = |A| > \ap_0$.  Hence
there exist $C \subseteq A\bs A'$, and $n \in \N$ such that $|C| =
|A\bs A'| = |A|$, and $|x(f(x))| \geq 1/n$ for all $x \in C$.  Now $|C|
= |A| > |\G'| \geq |f(C)|$.  Therefore, there is a $\g_0 \in f(C)$ such
that $D = f^{-1}\{\g_0\} \cap C$ is infinite.  Note that $x \in D$ implies
$|x(\g_0)| \geq 1/n$.  Now for any finite subset $F$ of $D$, 
\[ \|\sum_{x\in F}\sgn x(\g_0)x\| \geq \sum_{x\in
F}|x(\g_0)|\|\chi_{\{\g_0\}}\| \geq \frac{|F|}{n}w(\g_0)^{\frac{1}{p}} . \]
As $D$ is infinite, this violates condition (\ref{dom}).

Now for each $x \in \wlp(\G,w)$, let $\supp x = \{\g \in \G : x(\g) \neq
0\}$. Clearly $|\supp x| \leq \aleph_0$.  Therefore, $|\bc_{x\in B}\supp
x| \leq \max\{|B|,\aleph_0\} < |A|$.   Let $\G_1 = \bc_{x\in B}\supp x$.
By the above, there is a subset $A_1$ of $A$, having
the 
same cardinality as $A$, such that $x\chi_{\G_1} = 0$ for all $x \in
A_1$. 
It remains to choose a pairwise
disjoint subset of $A_1$ of cardinality $|A|$.  This will be done by
induction. Choose $x_0$ arbitrarily in $A_1$.  Now suppose a pairwise
disjoint collection
$(x_\rho)_{\rho<\beta}$ has been chosen up to some ordinal $\beta <
|A| = |A_1|$. Since $|A_1|$ is a cardinal, $|\beta| < |A_1|$.  Hence
$|\bc_{\rho<\beta}\supp x_\rho| \leq \max\{|\beta|,\aleph_0\} < |A_1|$.
Let $\G_2 = \bc_{\rho<\beta}\supp x_\rho$. Using the first part of the
proof again, we find a $x_\beta \in A_1$ such that
$x_\beta\chi_{\G_2} = 0$.  It is clear that the collection
$(x_\rho)_{\rho\leq\beta}$ is pairwise disjoint.  This completes the
inductive argument.  Consequently, we obtain a pairwise disjoint
collection $C = (x_\rho)_{\rho<|A|}$ in $A_1$. As each $x \in C$ is
disjoint from each $b \in B$, the proof is complete. 
\end{pf}

\begin{thm}\label{l2I}
Let $I$ and $\G$ be arbitrary sets such that $I$ is uncountable.
For any weight function $w$ on $\G$, and any $p \neq 2$, 
$1 < p < \infty$, $\wlp(\G,w)$ does not contain
a subspace isomorphic to $\ell^2(I)$. 
Consequently, Theorem \ref{short} holds if $p \neq 2$. 
\end{thm}

\begin{pf}
For any set $I$, and any $i \in I$, let $\ep_i : \two^I \to \two$ be
the projection onto the $i$th coordinate. By Khinchine's inequality,
$(\ep_i)_{i\in I} \subseteq \wLp(\two^I)$ 
is equivalent to the unit vector basis of $\ell^2(I)$.
Hence the first statement of the theorem implies the second.
Now suppose $(x_i)_{i\in I}$ is a set of normalized elements of
$\wlp(\G,w)$ which is equivalent to 
the unit vector basis of $\ell^2(I)$. If $I$ is uncountable, apply 
Proposition \ref{disj} with $A = I$, $B = \emptyset$ 
to obtain an uncountable
$C \subseteq I$ such that $(x_i)_{i\in C}$ are pairwise disjoint. 
Since $\wlp(\G,w)$ satisfies an upper $p$-estimate, there is a constant
$0 < K < \infty$ such that 
\[ K^{-1}|F|^{\frac{1}{2}} \leq \bl\|\sum_{i\in F}x_i\br\| \leq
K|F|^{\frac{1}{p}} \]
for every finite subset $F$ of $C$. We conclude that $1 < p < 2$.
Denote by $\mu$ the measure associated with $(\G,w)$. 
For each $i \in C$,
there is a rational number $c_i > 0$ such that 
$c_i(\mu\{|x_i|>c_i\})^{\frac{1}{p}} > \frac{1}{2}$.
By using an uncountable subset of $C$ if necessary, we may assume that 
$c_i = c$, a constant, for all $i \in C$. For any finite subset $F$ 
of $C$,
\[ \mu\bl\{\bl|\sum_{i\in F}x_i\br|>c\br\} = 
\sum_{i\in F}\mu\{|x_i|>c\} > (2c)^{-p}|F| . \]
Hence $\|\sum_{i\in F}x_i\| > \frac{1}{2}|F|^{\frac{1}{p}}$. 
Since $1 < p < 2$, and $(x_i)_{i\in C}$ is equivalent to the unit
vector basis of $\ell^2(C)$, we have reached a contradiction.
\end{pf}

The proof of Theorem \ref{short} for the case $p = 2$ is more involved.
Let $(h_n)$ denote the $L^\infty$-normalized Haar functions on $[0,1]$
\cite[Definition 1.a.4]{LT}. Then by \cite[Theorem 2.c.6]{LT2},
$(h_n)$ is an unconditional basis of $\Mp[0,1]$. We first show that if
$T : \Mt[0,1] \to \wLt\msp$ is an embedding, then $(Th_n)$ cannot be
pairwise disjoint.

\begin{prop}\label{Haar}
Suppose $T : \Mt[0,1] \to \wLt\msp$
is an embedding for some measure space $\msp$.
Then $(Th_n)$ cannot be a pairwise disjoint sequence.
\end{prop}

\begin{pf}
By Theorem 1.d.6(ii) in \cite{LT2}, there is a constant $D$ such that
\[ \|\sum a_nh_n\| \geq D^{-1}\|(\sum|a_nh_n|^2)^{1/2}\| \]
for every
sequence of scalars $(a_n)$ which is eventually zero. Given $m \in \N$,
let $a_{ij} = i+j-1 \pmod{2^m}$, $1 \leq i, j \leq 2^m$. For $1 \leq j
\leq 2^m$, define
\[ g_j = \sum^{2^m}_{i=1}\bl(\frac{2^m}{a_{ij}}\br)^{1/2}\chi_{I_i} , \]
where $I_{i} = [\frac{i-1}{2^m},\frac{i}{2^m})$. 
If $1 \leq j \leq 2^m$, there exists 
\[ f_j = \sum^{2^{m+j}}_{n=2^{m+j-1}+1}\!\!b_nh_n \in 
\spn\{h_n : 2^{m+j-1} < n \leq 2^{m+j}\}\]
such that $|f_j| = g_j$. Note that $(f_j)^{2^m}_{j=1}$
is a normalized sequence in $\Mt[0,1]$.  
If $T : \Mt[0,1] \to \wLt\msp$ is a bounded
linear operator such that $(Th_n)$ is pairwise disjoint, then
$(Tf_j)^{2^m}_{j=1}$ is a pairwise disjoint sequence which is bounded
in norm by $\|T\|$. Hence, using the upper $2$-estimate in $\wLt\msp$,
$\|\sum^{2^m}_{j=1}Tf_j\| \leq 2^{m/2}\|T\|$.  On the other hand,
\begin{align*}
\bl\|\sum^{2^m}_{j=1}f_j\br\| &=
 \bl\|\sum^{2^m}_{j=1}\sum^{2^{m+j}}_{n=2^{m+j-1}+1}\!\!b_nh_n\br\| \\
&\geq
 D^{-1}\bl\|\bl(\sum^{2^m}_{j=1}
 \sum^{2^{m+j}}_{n=2^{m+j-1}+1}\!\!|b_nh_n|^2\br)^{1/2}\br\|
 \\
&= D^{-1}\bl\|\bl(\sum^{2^m}_{j=1}|f_j|^2\br)^{1/2}\br\| \\
&= D^{-1}\bl\|\bl(\sum^{2^m}_{j=1}|g_j|^2\br)^{1/2}\br\| \\
&= D^{-1}\bl(2^m\sum^{2^m}_{j=1}\frac{1}{j}\br)^{1/2}\|\chi_{[0,1)}\| 
= D^{-1}\bl(2^m\sum^{2^m}_{j=1}\frac{1}{j}\br)^{1/2} .
\end{align*}
Since $m$ is arbitrary, $T$ cannot be an embedding.
\end{pf}

We now complete the proof of Theorem \ref{short} for the case $p = 2$.
Suppose $I$ is uncountable, and $T : \wLt(\two^I) \to \wlt(\G,w)$  
is a bounded linear operator. We will construct a
sequence $(g_n) \subseteq \wLt(\two^I)$ 
which is equivalent to the Haar basis $(h_n) \subseteq
\wLt[0,1]$, and such that $(Tg_n)$ is a pairwise disjoint sequence in
$\wlt(\G,w)$.  An appeal to Proposition \ref{Haar} will then yield the
desired result that $T$ is not an embedding. 
Let the functions $(\ep_i) \subseteq \wLt(\two^I)$ be the same as those
which appeared in the proof of Theorem \ref{l2I}. For any finite subset
$F$ of $I$, and any $\two$-valued sequence $(b_i)_{i\in F}$, the family
\[ \bl(\bl(\prod_{i\in F}\chi_{\{\ep_i=b_i\}}\br)
\ep_j\br)_{j\in I\backslash F} \]
is easily seen to be equivalent to the unit vector basis of
$\ell^2(I\backslash F)$.  Applying
Proposition \ref{disj}, we see that for any $F$ and $(b_i)_{i\in F}$ as
above, and any countable $S \subseteq \G$, there exists $j \in
I\backslash F$ such that $(T((\prod_{i\in
F}\chi_{\{\ep_i=b_i\}})\ep_j)){\chi_S} = 0$.  For any $x \in
\wlt(\G,w)$, we let its {\em support}\/ be the set $\supp x = \{\g\in\G
: x(\g) \neq 0\}$. Any element in $\wlt(\G,w)$ has countable support. 
Let $g_1$ be the
identically $1$ function on $\two^I$. Then there exists $i_2 \in I$
such that $\supp T\ep_{i_2} \cap\, \supp Tg_1 = \emptyset$. Let $g_2 =
\ep_{i_2}$. Now $\supp Tg_1 \cup\, \supp Tg_2$ is countable. Therefore,
one can find $i_3 \neq i_2$ such that $\supp
T((\chi_{\{\ep_{i_2}=-1\}})\ep_{i_3})$ is disjoint from $\supp Tg_1
\cup \supp Tg_2$.  Define $g_3 = (\chi_{\{\ep_{i_2}=-1\}})\ep_{i_3}$.
Next define $g_4 = (\chi_{\{\ep_{i_2}=1\}})\ep_{i_4}$, where $i_4$ is
chosen so that it is distinct from $i_2, i_3$, and $\supp Tg_4$ is
disjoint from $\cup^3_{n=1}\supp Tg_n$. Continuing in this way, we
obtain the desired sequence $(g_n)$.

\section{The space $\wLp(\dsum J_\al)$}

In this section, we present the proof of Theorem \ref{long}. Let $A$
be an uncountable set, and let $w$ be
a weight function defined on a set $\G$. 
Suppose $T : \Mp(\dsum J_\al) \to
\wlp(\G,w)$ is a bounded linear operator. The first step is to show
that the range of $T$ is mostly contained in $\mpi(\G,w)$. This will
require the following technical lemma.

\begin{lem}\label{tech}
Let $k \in \N$ be given, and let $\delta, c_1, c_2, \dots, c_k$ be
strictly positive numbers. Suppose $l \in \N$ is so large that  
\[ l\bl(\frac{\min(c_1,\dots,c_k)}{\max(c_1,\dots,c_k)}\br)^p \geq 1 . \]
Let $\msp$ be any measure space, and 
let $f_1, \dots, f_l$ be pairwise disjoint functions in $\wLp\msp$
such that $|f_j| \geq \sum^k_{m=1}c_m\chi_{\si(m,j)}$, where, for each
$j$, $\si(1,j),\dots,\si(k,j)$ are pairwise disjoint sets in $\Si$ such
that $\mu(\si(m,j)) > (\delta/c_m)^p,\, 1 \leq m \leq k$.  Then
\[ \bl\|\sum^l_{j=1}j^{-1/p}f_j\br\| \geq 
\bl(\frac{k}{2}\br)^{1/p}\delta . \]
\end{lem}

\begin{pf}
We may assume without loss of generality that 
$c_1 \geq \dots \geq c_k > 0$.
Then $l(c_k/c_1)^p \geq 1$.  For $1 \leq m \leq k$, let $i_m$ be the
largest integer in $\N$ $\leq l(c_m/c_1)^p$. Note that 
$1 \leq i_m \leq l$ .
For any $\ep < c_1l^{-1/p}$, 
\begin{align*}
\bl\{\bl|\sum^l_{j=1}j^{-1/p}f_j\br| > \ep\br\} &\supseteq
 \bigcup\{\si(m,j) : c_mj^{-1/p} > \ep\} \\
&\supseteq \bigcup\{\si(m,j) : c_mj^{-1/p} \geq c_1l^{-1/p}\} \\
&= \bigcup\{\si(m,j) : j \leq l(c_m/c_1)^p\} \\
&= \bigcup^k_{m=1}\bigcup^{i_m}_{j=1}\si(m,j) .
\end{align*}
Thus
\begin{align*}
\mu\bl\{\bl|\sum^l_{j=1}j^{-1/p}f_j\br| > \ep\br\} &\geq  
 \sum^k_{m=1}\sum^{i_m}_{j=1}\mu(\si(m,j)) \\
&> \sum^k_{m=1}\sum^{i_m}_{j=1}(\delta/c_m)^p = 
 \sum^k_{m=1}i_m(\delta/c_m)^p .
\end{align*}
Now $i_m \geq 1$ implies $i_m \geq (1+i_m)/2 \geq 2^{-1}l(c_m/c_1)^p$.
Hence 
\[
\mu\bl\{\bl|\sum^l_{j=1}j^{-1/p}f_j\br| > \ep\br\} \geq
 \sum^k_{m=1}2^{-1}l(c_m/c_1)^p(\delta/c_m)^p 
= (lk/2)(\delta/c_1)^p . \]
Therefore, 
\[ \bl\|\sum^l_{j=1}j^{-1/p}f_j\br\| \geq 
\ep\bl(\mu\bl\{\bl|\sum^l_{j=1}j^{-1/p}f_j\br| > \ep\br\}\br)^{1/p} 
\geq (\ep\delta/c_1)(lk/2)^{1/p} . \]
Taking the supremum over all $\ep < c_1l^{-1/p}$ yields the desired result.
\end{pf}

\begin{prop}\label{range}
Let $A$ be an index set, and
let $T : \mpi(A) \to \wlp(\G,w)$ be a bounded linear operator 
for some $(\G,w)$.
Then $T\chi_{\{\al\}} \in \mpi(\G,w)$ for all but countably many 
$\al \in A$.
\end{prop}   

\begin{pf}
Let $f_\al = T\chi_{\{\al\}}$, and assume $f_\al \notin \mpi(\G,w)$
for uncountably many $\al$.  
Applying Proposition \ref{disj}, we may assume that the $f_\al$'s are 
pairwise disjoint. 
Choose an uncountable $A_0 \subseteq A$, and 
$\delta > 0$, such that $d(f_\al,\mpi(\G,w)) > \delta$ for all $\al \in
A_0$. For each $\al \in A_0$, there is a rational $r > 0$ 
such that $\mu\{|f_\al| > r\} > (\delta/r)^p$, where $\mu$ is the measure 
associated with $(\G,w)$. Hence we can find an uncountable 
$A_1 \subseteq A_0$, and $c_1 > 0$ such that 
$\mu\{|f_\al| > c_1\} > (\delta/c_1)^p$ for every $\al \in A_1$. 
For all $\al \in A_1$, choose a finite set $\si(1,\al) \subseteq
\{|f_\al| > c_1\}$ such that $\mu(\si(1,\al)) > (\delta/c_1)^p$.
Now $\|f_\al - f_\al\chi_{\si(1,\al)}\| > \delta$ for all 
$\al \in A_1$. Arguing as before, we find an uncountable $A_2 \subseteq
A_1$, and $c_2 > 0$ such that 
\[ \mu\{|f_\al - f_\al\chi_{\si(1,\al)}| > c_2\} > (\delta/c_2)^p \]
for all $\al \in A_2$. Hence, for each $\al\in A_2$, 
there exists a finite set $\si(2,\al) \subseteq
\{|f_\al| > c_2\}$, disjoint from $\si(1,\al)$,
such that $\mu(\si(2,\al)) > (\delta/c_2)^p$.
Continue inductively to obtain a decreasing sequence of uncountable subsets 
$(A_m)$ of $A$, a positive sequence $(c_m)$, and finite subsets
$\si(m,\al) \subseteq \{|f_\al| > c_m\}$ 
for all $\al \in A_m$, such that $\mu(\si(m,\al)) > (\delta/c_m)^p$, and
$\si(m,\al)\cap\si(n,\al) = \emptyset$ if 
$\al \in A_m\cap A_n$ and $m \neq n$.
Now let $k \in \N$ be given. Choose $l$ so large that 
\[ l\bl(\frac{\min(c_1,\dots,c_k)}{\max(c_1,\dots,c_k)}\br)^p \geq 1 . \]
Lemma \ref{tech} 
implies that $\|\sum^l_{j=1}j^{-1/p}f_j\| \geq (k/2)^{1/p}\delta$ if
$\al_1,\dots,\al_l$ are distinct elements of $A_k$. 
This violates the boundedness of $T$ since $k$ is arbitrary.
\end{pf}

\begin{cor}\label{range2}
Let $A$ be an index set. For each $\al \in A$, $n \in \N$, and $1 \leq j
\leq 2^n$, let $f_{n,j,\al}$ be the characteristic function of the 
subinterval $[\frac{j-1}{2^n},\frac{j}{2^n})$ in $J_\al$. If $T : 
\Mp(\dsum J_\al) \to \wlp(\G,w)$ is a bounded linear operator, then 
all but countably many members of $\{Tf_{n,j,\al} : \al\in A, n \in \N,
1 \leq j \leq 2^n\}$ belong to $\mpi(\G,w)$.
\end{cor}

\begin{pf}
If $n$ is fixed, the collection $\{f_{n,j,\al} : \al\in A,\, 
1 \leq j \leq 2^n\}$ is equivalent to the unit vector basis in 
$\mpi(A\times\{1,\dots,2^n\})$. Apply Proposition \ref{range} to
complete the proof.
\end{pf}

If $A$ is uncountable, and $T : \wLp(\dsum J_\al) \to \wlp(\G,w)$
is an embedding, then it follows from Corollary \ref{range2} that there
exists $\al_0 \in A$ such that $T(\Mp(J_{\al_0})) \subseteq \mpi(\G,w)$, 
where we identify $\Mp(J_{\al_0})$ with a subspace of $\wLp(\dsum J_\al)$
in the obvious way.  Hence $\Mp[0,1]$ embeds into $\mpi(\G,w)$.
The proof of Theorem \ref{long} is completed by showing that this is 
impossible.  Once again, we find it necessary to distinguish between the 
cases $p \neq 2$ and $p = 2$. If $p \neq 2$, we use a 
Kadec-Pe\l czy\'{n}ski type argument \cite{KP} to show that $\ell^2$ does
not embed into $\mpi(\G,w)$. For $p = 2$, we resort once again to 
Proposition
\ref{Haar}. If $f$ is a real valued function and $1 <M < \infty$, let 
$(f)_M = f\chi_{\{M^{-1} < |f| < M\}}$.

\begin{lem}\label{sum}
Let $\msp$ be any measure space, and suppose $1 < p < \infty$. If $(f_n)$ 
is a pairwise disjoint sequence in the unit ball of 
$\wLp\msp$, and $(M_n)$ is a real 
sequence such that $1 < M_n \leq 2^{-1/p}M_{n+1}$ for all $n \in \N$,
define $g_1 = (f_1)_{M_1}$ and 
\[ g_{n+1} = (f_{n+1})_{M_{n+1}}-(f_{n+1})_{M_n} . \]
Then $ \sup_k\bl\|\sum^k_{n=1}g_n\br\| \leq 4$.
\end{lem}

\begin{pf}
Let $g = \mbox{pointwise-}\!\sum g_{n}$ and 
$M_0 = 1/{M_1}$. If $M_{k-1} \leq c < M_k$ for some $k \in \N$, then
\begin{align*}
\mu\{|g| > c\} &= \sum^\infty_{n=k}\mu\{|g_n| > c\} \\
 &= \mu\{|g_k| > c\} + \sum^\infty_{n=k+1}\mu\{|g_n| \geq M_{n-1}\} \\
 &\leq c^{-p} + \sum^\infty_{n=k+1}{M^{-p}_{n-1}} 
  \leq c^{-p} + 2{M^{-p}_k} \leq 3c^{-p} .
\end{align*}
On the other hand, if $M^{-1}_{k+1} \leq c < M^{-1}_k$ for some 
$k \in \N$, then
\begin{align*}
\mu\{|g| > c\} &= \sum^k_{n=1}\mu\{|g_n| > M^{-1}_n\} 
   + \mu\{|g_{k+1}| > c\} + \sum^\infty_{n=k+2}\mu\{|g_n| \geq M_{n-1}\} \\
 &\leq \sum^k_{n=1}M^p_n + c^{-p} + \sum^\infty_{n=k+2}M^{-p}_{n-1} \\
 &\leq 2M^p_k + c^{-p} + M^{-p}_{k} \leq 2c^{-p} + c^{-p} + 1
  \leq 4c^{-p} .
\end{align*}
Hence $g \in \wLp\msp$, and $\|g\| \leq 4$. 
\end{pf}

\begin{thm}\label{ell2}
For any $(\G,w)$, and $1 < p < \infty$, $p \neq 2$, there is no embedding
of $\ell^2$ into $\mpi(\G,w)$. 
\end{thm}

\begin{pf}
Suppose, on the contrary, that $\mpi(\G,w)$ contains a sequence 
$(f_n)$ equivalent
to the unit vector basis of $\ell^2$. Since each $f_n$ has countable 
support, we may assume that $\G$ is countable. Then clearly 
$(\chi_{\{\g\}})_{\g\in\G}$ is an unconditional basis of $\mpi(\G,w)$.
Since $(f_n)$ is a weakly null sequence, we may apply the 
Bessaga-Pe\l czy\'{n}ski 
selection principle \cite[Proposition 1.a.12]{LT} to it. 
Thus, we may assume without loss of generality that $(f_n)$ is pairwise 
disjoint. Since $\mpi(\G,w)$ satisfies an upper $p$-estimate, this is 
possible only if $1 < p < 2$. Now suppose there exists $1 < M < \infty$ 
such that $\limsup_n\|(f_n)_M\| > 0$. We may assume that there exists 
$\ep > 0$ such that $\|(f_n)_M\| > \ep$ for all $n$. For each $n$, choose 
$c_n \in [M^{-1},M]$ such that $c_n(\mu\{|f_n| > c_n\})^{1/p} > \ep$, where
$\mu$ is the measure associated with $(\G,w)$. Using the compactness of 
$[M^{-1},M]$, and going to a subsequence if necessary, we may assume the
existence of a $c \in [M^{-1},M]$ such that $c(\mu\{|f_n| > c\})^{1/p} > 
\ep$ for all $n$. Then
\[ \bl\|\sum^k_{n=1}f_n\br\| \geq 
 c\bl(\sum^k_{n=1}\mu\{|f_n| > c\}\br)^{1/p} \geq \ep k^{1/p} \]
for all $k \in \N$, a contradiction. Therefore, it must be that
$\lim_n\|(f_n)_M\| = 0$ for all $1 < M < \infty$. Note that 
$\lim_{M\to\infty}\!\|f - (f)_M\| = 0$ for every $f \in \mpi(\G,w)$.
By a standard perturbation argument, we obtain a subsequence of 
$(f_n)$, denoted again by $(f_n)$, and a real sequence $(M_n)$
satisfying  $1 < M_n \leq 2^{-1/p}M_{n+1}$ for all $n \in \N$, 
such that $(f_{n+1})$ is equivalent to 
$((f_{n+1})_{M_{n+1}}-(f_{n+1})_{M_n})$. 
Lemma \ref{sum}, however, shows that $(f_{n+1})$ cannot be equivalent to
the unit vector basis of $\ell^2$.
\end{pf}

We now give the proof of Theorem \ref{long}.
Assume that for some uncountable set $A$, 
$\wLp(\dsum J_\al)$ embeds into $\wlp(\G,w)$ for some $(\G,w)$.
As in the 
discussion following Corollary \ref{range2}, $\Mp[0,1]$ embeds into 
$\mpi(\G,w)$. Since the sequence of Rademacher functions in 
$\Mp[0,1]$ is equivalent to the unit vector basis of $\ell^2$, Theorem 
\ref{ell2} implies that this is impossible unless $p = 2$.  Now let 
$T : \Mt[0,1] \to \mti(\G,w)$ be an embedding. Without loss
of generality, assume that $\|Tf\| \geq \|f\|$ for all $f \in \Mt[0,1]$.
Denote by $(r_n)$, respectively $(h_n)$, 
the sequence of Rademacher functions, respectively Haar functions, 
on $[0,1]$. Note that for all $f \in \Mt[0,1]$, 
$f\cdot r_n \to 0$ weakly as $n \to \infty$.
Let $f_1 = |h_1|$.
If $k \in \N$, and $2^{k-1} < j \leq 2^k$, let $f_j = \sqrt{2^{k-1}}|h_j|$. 
Define $n_1 = 1$. Since $x_1 = T(h_1\cdot r_{n_1}) \in \mti(\G,w)$, 
there is 
a finite subset $\si_1$ of $\G$ such that $\|x_1\chi_{\si^c_1}\| 
< 2^{-3}$. Now suppose that numbers $n_i$ and finite sets 
$\si_i$ have been chosen for $i \leq j$. 
Since $T(f_{j+1}\cdot r_n) \to 0$ weakly as $n \to \infty$, and 
$\cup^j_{i=1}\si_i$ is finite, there exists $n_{j+1} > n_j$ so that 
$\|x_{j+1}\chi_{\cup^j_{i=1}\si_i}\| < {2^{-j-4}}$, where $x_{j+1}
= T(f_{j+1}\cdot r_{n_{j+1}})$. Now we can choose a finite subset
$\si_{j+1}$ of $\G$, disjoint from $\cup^j_{i=1}\si_i$, such that 
$\|x_{j+1}\chi_{\si^c_{j+1}}\| < {2^{-j-3}}$. 
Finally, let $y_j = x_j\chi_{\si_j}$ for all $j \in \N$. 
Then $(y_j)$ is pairwise disjoint sequence, and hence is a basic 
sequence with
basis constant $1$. Moreover, 
\[ \|y_j\| >  \|x_j\| - 2^{-j-2} \geq 
\|f_j\cdot r_{n_j}\| - 2^{-j-2} > 1/2 . \]
Also, $\sum\|x_j - y_j\| < 1/4$. By Proposition 1.a.9
in \cite{LT}, $(y_j)$ and $(x_j)$ are equivalent.
But then $(f_j\cdot r_{n_{j}})$ is equivalent to a pairwise disjoint 
sequence in $\wlp(\G,w)$. However, it is easy to see that 
$(f_j\cdot r_{n_{j}})$ is equivalent to $(a_jh_j)$, where 
$a_1 = 1$ and 
$a_j = \sqrt{2^{k-1}}$ if $2^{k-1} < j \leq 2^k$. Hence we obtain an 
embedding $S$ of $[(h_j)]$ into $\wlp(\G,w)$ such that $(Sh_j)$ 
is a pairwise 
disjoint sequence. As $(h_j)$ is a basis of $\Mt[0,1]$, we have reached
a contradiction to Proposition \ref{Haar}.  This completes the proof
of Theorem \ref{long}.


\end{document}